\documentclass[a4paper, 11pt]{article}

\usepackage[utf8]{inputenc}
\usepackage{amsmath}
\usepackage{amsthm}
\usepackage{amsfonts}
\usepackage{amssymb}
\usepackage{amstext}
\usepackage{dsfont}
\usepackage{graphicx}
\usepackage{color}
\usepackage[colorlinks]{hyperref}
\usepackage{epigraph}
\usepackage[left=3.5cm,top=2cm,bottom=3.5cm,right=3cm,nohead,nofoot]{geometry}
\usepackage{tikz}
\usepackage{tkz-graph}
\usepackage{tkz-berge}
\usetikzlibrary{arrows,shapes,matrix}
\graphicspath{{Figures/}{.}}
\allowdisplaybreaks[4]
\usepackage{enumitem}

\title{Asymptotic properties of the density of particles\\ in $\beta$-ensembles}
\author{M. Dal Borgo, E. Hovhannisyan, A. Rouault} 
\theoremstyle{plain}

\newtheorem{thm}{Theorem}[section]
\newtheorem{proposition}[thm]{Proposition}

\newtheorem{conjecture}[thm]{Conjecture}
\newtheorem{definition}[thm]{Definition}

\newtheorem{rmk}[thm]{Remark}

\numberwithin{equation}{section}
\def\half{\frac{1}{2}}
\def\betah{\frac{\beta}{2}}
\def\l{\lambda}

\def\m{\mu}

\def\b{\beta}
\def\d{\delta}
\def\g{\gamma}
\def\G{\Gamma}

\def\dd{\mathrm {d}}
\def\mm{\mathrm {m}}

\def\N{{\mathbb N}}

\def\R{{\mathbb R}}
\def\C{{\mathbb C}}

\def\P{{\mathbb P}}
\def\E{{\mathbb E}}

\newcommand{\Int}{\operatorname{int}\,}
\newcommand{\Clo}{\operatorname{clo}\,}

\begin{document}
\maketitle
\begin{abstract}
We extend recent results on the Asymptotic Equipartition Property for the density of $n$ particles in $\b$-ensembles, as $n$ tends to infinity. We prove the Large Deviation Principle of the log-density for a general potential and the mod-gaussian convergence in the classical examples.
\end{abstract}

\bigskip

\section{Introduction}
We study the asymptotic properties of the $\log$-density in the $\b$-ensembles on $\mathbb R$ or on the unit circle $\mathbb T$, in the regime of large number of particles. We consider a system of $n$ random particles with a density of the form
\begin{align}\label{V_density}
{\bf P}_{n, \beta}^V\left(  \lambda_1, \dots, \lambda_n\right) = \frac{1}{\mathcal{Z}_n^V(\beta)} \prod_{1\le j<k\le n} \left| \lambda_k - \lambda_j \right|^\beta \prod_{k=1}^n e^{- \frac{\beta}{2}  n V(\lambda_k)},
\end{align}
with respect to the Lebesgue measure $d\lambda_1\dots d\lambda_n$. Here
 $\b>0$ is the inverse temperature, the function $V$ represents the potential and $\mathcal{Z}_n^V(\beta)$ is the normalization factor  (partition function) given by
\begin{align}
\label{partV}{\mathcal Z}_n^V(\beta) = \int_{\mathbb R^n}   \prod_{1\le j<k\le n} |\lambda_j - \lambda_k
|^\beta  e^{- n \frac{\beta}{2} \sum_{k=1}^n V(\lambda_k)}\ \   d\lambda_1 \dots d\lambda_n\,.\end{align}
This model comes from random matrix theory where $\lambda_1, \dots, \lambda_n$ are the eigenvalues of a random matrix, orthogonal, Hermitian, or  self-dual according to $\beta= 1,2,4$ whose distribution is invariant by unitary transformation.

We restrict ourselves to the case where the potential $V$ is a real analytic function satisfying the condition
\begin{align} \label{confin}
V(\lambda) \geq 2(1+\varepsilon) \log \left( 1+|\lambda|\right),
\end{align}
for all sufficiently large enough $\lambda$. We should stress that 
the extra factor $n$ in the exponential term of (\ref{V_density}) ensures that both the confining  potential 
$n\sum_{k=1}^n V(\l_k)$ and the repulsive potential term $\sum_{1\le j<k\le n}\left| \lambda_k - \lambda_j \right|$ are of the same order.

On $\mathbb T$, we have the same formulas if we replace $\lambda_k$ by $e^{i\theta_k}$ and $\mathbb R^n$ by $\mathbb T^n$ and if we assume that $V$ is continuous on $\mathbb T \setminus \{1\}$ and semicontinuous at $1$.

A huge literature deals with the asymptotic behavior of the  empirical distribution, the extreme values, the spacing between particles , etc.

Here we are concerned with the density (\ref{V_density}) as a random object, just like the likelihood in statistics. 

Our results are mostly based on the article of Bufetov et al. \cite{BMSS13}, where  
the asymptotic equipartition property (AEP) have been obtained for  some ensembles. The AEP 
 states that the logarithm of the density function, 
\[\mathcal L_{n, \beta}^V (\lambda) :=  \log {\bf P}_{n, \beta}^V (\lambda)\]
converges almost surely, after some normalization, to a constant which can be seen as minus the entropy of the corresponding measure, (see Section 2.1.1 in \cite{BMSS13} to the connection with the Shannon-McMillan-Breiman theorem). 

Moreover Bufetov et al. proved a central limit theorem (CLT). In general, the AEP is strongly connected to the existence of a large deviation principle  (LDP) and in Gibbs models, to a suitable expansion of the partition function. 
After the large amount of studies initiated by the seminal paper \cite{AG97}, Serfaty and Lebl\'e \cite{LS15} proved an LDP for a microscopic rescaled empirical function. On the way, the authors conjectured that a  unicity property could lead directly to the results of \cite{BMSS13}, with weaker assumptions.

The goal of this paper,  is to give fine asymptotic results of the log-density such as  large deviations, Berry-Esseen estimates, local and extended CLT. To obtain some of these results we use a framework of mod-Gaussian convergence that has been created and developed in~\cite{ FMN16, KN12, JKN11, DKN15}. Also large deviations at the usual scale are derived using the G\"artner-Ellis theorem. 

It is worth mentioning that in \cite{POP07}, Popescu studied  the empirical energy of the  $\b$-Hermite ensemble  and proved the CLT  as well as the LDP.

So, our results have to be seen as a refinement of   the results of \cite{BMSS13} 
and of \cite{POP07}.

The structure of the paper goes as follows. 
In Section \ref{basicex} we present the examples which are particularly suitable for our study and we explain in Section \ref{IS} why they are convenient. In Section \ref{section_2}
 we recall the definition of LDP and define mod-Gaussian convergence and the limiting theorems this convergence implies. In Section \ref{subs}, we introduce our results the eigenvalue density 
of the ensembles. Section \ref{proofdens} contains the proofs related to Section \ref{subs}.  An appendix (Section \ref{appendix}) gathers the properties of the special functions used in this paper.
\section{Basic examples}
\label{basicex}
We present here the basic examples which allow explicit computations of partition functions.
\subsection{Real line}
1) The Hermite (or Gaussian) ensemble corresponds to 
\begin{align}
\label{defpotH}V(x) := H(x) := \frac{x^2}{2}.
\end{align}
Its equilibrium measure $\mu_H$ is the semicircle law 
with the support ${\bf S} = [-2, 2]$ and the density is
\begin{equation}
\label{defSC}\rho_H (x) = \frac{1}{2\pi} \sqrt{4-x^2}\,.\end{equation}

\noindent2)  The Laguerre ensemble of parameter $\theta \geq 1$ corresponds to  
\begin{align}
\label{defpotL}V(x) := L(x) := \theta x - (\theta -1) \log x \ , \ ( x > 0)\,,
\end{align}
(see \cite{Mehta} formula 17.6.5 or  \cite{F10} Proposition 4.7.3).\footnote{Notice that there are various scalings in the literature.}
 
\noindent Its equilibrium measure $\mu_L$ is the  Marchenko-Pastur distribution with the support ${\bf S} = [c_-, c_+]$ where 
\[c_\pm = \frac{(1\pm \sqrt \theta)^2}{\theta}\,,\] and the density is 
\begin{equation}
\label{defMP}\rho_L(x)  = \theta\frac{\sqrt{(c_+ -x)(x-c_-)}}{2\pi x}\,.\end{equation}

\noindent 3) The Jacobi ensemble of parameters $\kappa_1, \kappa_2 > 0$ corresponds to 
\begin{align}
\label{defpotJ}
V(x) &:= -\kappa_1 \log x - \kappa_2 \log (1-x) \ \ (x \in [0,1))
\,.
\end{align}
(see \cite{Mehta} formula 17.1.3 or  \cite{F10} formula (4.3)).

\noindent Its equilibrium measure is the Kesten-McKay distribution whose support is ${\bf S} = [a^-, a^+]$  where 
$$
a^\pm := \frac{1}{2} + \frac{\kappa_1^2 - \kappa_2^2 \pm 4 \sqrt{(1+\kappa_1)(1+\kappa_2)(1+\kappa_1+\kappa_2)}}{2(2+\kappa_1+\kappa_2)^2}\, 
$$and the density is 
\begin{equation}
\label{defJ}
\rho_{KMK}(x) = \frac{C_{a_-, a_+}}{2\pi}\frac{\sqrt{(a^+ -x)(x- a^-)}}{x(1-x)} \,.\end{equation}
\begin{equation}
\label{defca}
C_{a_-, a_+} = 2 + \kappa_1 + \kappa_2\,,
\end{equation} (see \cite{F10} (3.77)).

\noindent 4) The generalized Cauchy ensemble of parameter $\dd$ (see \cite{F10} p. 174) corresponds to 
\begin{align}
V(x) = (1+ \dd) \log (1+x^2)\,.
 \end{align}
\noindent Its equilibrium measure has a support  ${\bf S} = [-\mm, \mm]$ with $\mm^2  = (1 + 2\dd) \dd^{-2}$
and the density is
\begin{equation}
\label{defCy}
\rho_{Cy}(x) 
= \frac{1}{\pi(\sqrt{1+\mm^2} -1)}\frac{\sqrt{\mm^2-x^2}}{1+x^2}\,.
\end{equation}
\subsection{Unit circle}

1) The circular ensemble with $V =0$ 
whose equilibrium measure is the normalized Lebesgue measure on $\mathbb T$.

\noindent 2) The circular Jacobi ensemble (see \cite{F10} p.126),  also called Hua-Pickrell ensemble corresponds to
\begin{align}
\label{defpotHP}V(z)&= -2\dd \log |1-z|
\,.\end{align}
 In the case $\beta = 2$ the stereographic projection (or Cayley transform) realizes  a perfect bijection between the generalized Cauchy ensemble and the circular Jacobi ensemble.

\noindent Its equilibrium measure is supported by an arc $[\theta_\dd , 2\pi - \theta_\dd]$ 
where $
\theta_\dd \in (0, \pi)$ is such that 
\begin{equation}\sin \frac{\theta_\dd}{2}  = \frac{\dd}
{ 1 + \dd}\,, 
\end{equation}
and its density is 
 \begin{equation}
\label{limmeas}
\rho_{HP} (x)  = (1+\dd) \frac{\sqrt{\sin^2(\frac{\theta}{2}) -
\sin^2(\frac{\theta_{\dd}}{2})}}{2\pi \!\ \sin(\frac{\theta}{2})}\,,
\end{equation} 
(see \cite{BNR09}).

\section{Selberg formulas}
\label{IS}
Most of the results of this paper are obtained for models where the partition function is a product and quotient of Gamma functions. This is a remarkable consequence of the importance of the Selberg formula. For a complete treatment of this topic, we refer to 
 Chapter 4 of \cite{F10} and also to the nice paper \cite{forrester2008importance}.

\noindent 1. In the Hermite case
\begin{align}
 \label{SelbergH}
\mathcal Z_n^H(\beta)
=
(2\pi)^{n/2} (n\beta')^{\frac{1}{2}((\beta' - 1) n -\beta' n^2)} \prod_{j=1}^n \frac{\Gamma(1+\beta' j)}{\Gamma(1+\beta')}\,.
\end{align}
\medskip

\noindent 2. In the Laguerre case,
\begin{align}
\label{SelbergL}
\mathcal Z_n^L(\beta) = \left(n\beta' \theta\right)^{-n(\theta n \beta' + 1 - \beta')} \prod_{j=1}^n \frac{\Gamma(1+ \beta' j)\Gamma (1 + \beta' (\theta n-n +j-1))}{\Gamma(1 + \beta')}
\end{align}
(see \cite{Mehta} formula 17.6.5 or  \cite{F10} Proposition 4.7.3).
\medskip

\noindent 3. In the Jacobi case,
\begin{align}
\mathcal Z_n^J (\beta) = \prod_{j=0}^{n-1} \frac{\Gamma\left((\kappa_1n +j)\beta' +1\right) \Gamma\left((\kappa_2n +j)\beta' +1\right) \Gamma\left(1+ (j+1)\beta'\right)}{\Gamma\left([(\kappa_1+ \kappa_2 +1)n +j-1]\beta' +1 \right)  \Gamma\left(1 + \beta'\right)}
\,.
\end{align}
(see \cite{Mehta} formula 17.1.3 or  \cite{F10} formula (4.3)).
\medskip

\noindent 4. In the generalized Cauchy case,
\begin{align}\mathcal Z_n^{Cy} = 2^{\beta'n(n-1) -2[(1+\dd)\beta' n -1]n} \pi^n \prod_{j=0}^{n-1}\frac{\Gamma\left(\beta'(j + 2\dd n + 2) -1\right)\Gamma\left(\beta'(j+1) +1\right)}{\Gamma\left(\beta'(j + \dd n + 1) \right)^2 \Gamma(1+\beta')}
 \end{align}
(see \cite{F10} formula (4.156)).
\medskip

\noindent 5. In the circular case
\begin{align}
\label{SelbergC}
\mathcal Z_n^C(\beta) = (2\pi)^n \frac{\Gamma (1 + n\b')}{\Gamma(1+\b')^n}\,,
\end{align}
(see \cite{Mehta} formula (17.7.1) or \cite{F10} Prop. 2.8.7).
\medskip

\noindent 6. In the circular Jacobi case
\begin{align}
\label{HP}
\mathcal Z_n^{HP}(\beta) = (2\pi)^n \prod_{j=0}^{n-1}\frac{\Gamma\left(\beta'(j+2\dd n) +1\right)\Gamma\left(\beta'(j+1)+1\right)}{\Gamma\left(\beta'(j+\dd n) +1\right)^2 \Gamma(1+ \beta')}\,.\end{align}
(see \cite{F10} formula (4.4)).

\section{Large deviations, mod-Gaussian convergence and consequences}\label{section_2}

We assume that the reader is familiar with theory of large deviations. We just recall the definition. For details, we refer to the classical book \cite{demboz98}. 
Let $(u_n)$
be a increasing positive sequence of real numbers with $\lim_{n \rightarrow \infty}u_n=\infty$.
\begin{definition}
\label{dldp}
We say that a sequence $(Q_{n})$ of probability measures on a measurable
Hausdorff space $G$ equipped with its Borel $\sigma$-field satisfies the  LDP with rate function $I$ at scale $u_n$ if:
\begin{itemize}
\item[i)] $I$ is lower semicontinuous, with values in
$\mathbb R^{+}\cup\{+\infty\}$.
\item[ii)] For any measurable set $A$ of $G$:
$$-I(\Int A)\leq
\liminf_{N\rightarrow\infty}\frac{1}{u_n}\log Q_n(A)\leq
\limsup_{N\rightarrow\infty}\frac{1}{u_n}\log Q_n(A)\leq
-I({\Clo A}),$$
where $I(A)=\inf_{\xi\in A}I(\xi)$ and $\Int A$ (resp. $\Clo A$) is the
interior (resp. the closure) of $A$.
\end{itemize}
We say that the rate function $I$ is good if its level set
$\{x\in G:\; I(x)\leq a\}$ is  compact for any $a\geq 0$.
More generally, a sequence of $G$-valued random variables 
is said to satisfy the
LDP if the sequence of their distributions satisfies the LDP.
\end{definition}
The notion of mod-Gaussian convergence is more recent and we will give definitions.

The framework of mod-Gaussian convergence for a sequence of random variables has been developed by Delbaen, F\'eray, Jacod, Kowalski, M\'eliot and Nikeghbali  in~\cite{FMN16, KN12, JKN11, DKN15, FMN17}.
Roughly speaking, given a sequence of random variables not converging in distribution, we look for a proper renormalization of their moment generating functions (and not of the sequence itself) which forces this convergence to happen.
If this normalization writes as the moment generating function of a Gaussian random variable with mean zero and variance blowing up to infinity, we speak of mod-Gaussian convergence for the sequence of random variables we started with.
Although mod-Gaussian convergence entails
much more information, for the aim of this paper, it will be an instrument to deduce asymptotic results, like a central limit theorem, large deviations etc.
In what follows we summarize some of the limiting results that mod-Gaussian convergence implies (see \cite{FMN16} for the results below).\\
\\
For $-\infty \leq c < 0 < d\leq \infty$, set
$$S_{(c,d)}=\{z\in\C, c < \Re (z) < d\}\ , \ S_d = S_{-d , \infty}\,.$$

\begin{definition}
\label{def:modGauss}  Let $\left(X_n\right)_{n\in\N}$ be a sequence of real valued random variables, and $\phi_n(z)=\E\left[e^{zX_n}\right]$ be their moment generating functions, which we assume to all exist in a strip $S_{c,d}$. 
 We assume that there exists an analytic function $\psi(z)$ not vanishing on the real part of $S_{(c,d)}$, such that locally uniformly on $S_{(c,d)}$, 
$$\lim_{n\to\infty}\E\left[e^{zX_n}\right]e^{-t_n\frac{z^2}{2}}:=\lim_{n\to\infty}\psi_n(z)=\psi(z),$$
where $\left(t_n\right)_{n\in\N}$ is some sequence going to infinity. We then say that 
 $\left(X_n\right)_{n\in\N}$ converges mod-Gaussian on $S_{(c,d)},$ with parameters $t_n$ and limiting function $\psi$.
\end{definition}
\vskip 10pt
In the remaining part of this section, $\left(X_n\right)_{n\in\N}$ is a sequence converging mod-Gaussian on $S_{(c,d)},$ with parameters $t_n$ and limiting function $\psi$. In particular, we are interested in the whole panorama of fluctuations that the mod-Gaussian framework implies for the sequence $\left(X_n\right)_{n\in\N}$. 
\\
\begin{thm}[Precise moderate deviations at scale $t_n$, Theorem 4.2.1 in \cite{FMN16}]\ \label{ldp}\\
For $x\in (0,d)$,
\begin{equation*}
\P\left[X_n\ge t_nx\right]=\frac{e^{-t_n\frac{x^2}{2}}}{x\sqrt{2\pi t_n}}\psi(x)\left(1+o(1)\right),
\end{equation*}
and for $x\in(c,0)$,
\begin{equation*}
\P\left[X_n\le t_nx\right]=\frac{e^{-t_n\frac{x^2}{2}}}{|x|\sqrt{2\pi t_n}}\psi(x)\left(1+o(1)\right).
\end{equation*}
\end{thm}
\vskip 10pt
\begin{thm}[Central limit theorem at  scale $o(t_n)$, Theorem 4.3.1 in \cite{FMN16}]\ \\
\label{clt}
For  $y=o\left(\sqrt{t_n}\right)$,
\begin{equation*}
\P\left[\frac{X_n}{\sqrt{t_n}}\ge y\right]=\P\left[
\mathcal{N}(0,1)\ge y\right]\left(1+o(1)\right)=\frac{e^{-\frac{y^2}{2}}}{y\sqrt{2\pi}}\left(1+o(1)\right).
\end{equation*}
\end{thm}

Note that Theorem \ref{clt} immediately implies a \emph{central limit theorem} for the rescaled sequence $\left(\frac{X_n}{\sqrt{t_n}}\right)_{n\in\N}$ (taking $y=O(1)$).
Under the additional assumption of mod-Gaussian convergence with a \emph{zone of control}, it is possible to deduce the speed of convergence of the above mentioned CLT and also to state a local limit theorem. Such notion has been introduced in \cite[Definition 4, Remark 2]{Speed} in the context of mod-stable convergence.

\begin{definition} 
Consider the following statements.
\begin{enumerate}[label=(\textbf{Z\arabic*}),ref=(Z\arabic*)]
\item \label{Z1} Fix $v \ge 1$, $w>0$ and $\g\in\R$. There exists a zone of convergence $[-Dt_n^\g,Dt_n^\g]$, $D>0$, such that for all $\xi$ in this zone, $\left|\psi_n(i\xi)-1\right|$ is smaller than
\begin{equation}
K_1|\xi|^v e^{K_2 |\xi|^w} \nonumber
\end{equation}
for some positive constants $K_1$ and $K_2$, that are independent of $n$.
\item \label{Z2} One has 
\begin{equation}
w \geq 2;\qquad -\half<\g\le\frac{1}{w-2};\qquad D\le \left(\frac{1}{4K_2}\right)^{\frac{1}{w-2}}. \nonumber
\end{equation}
\end{enumerate}
If Conditions (Z1) holds for some parameters $\g>-\half$ and $\nu, w, D, K_1, K_2$, then 
(Z2) can always be forced by increasing $w$, and then decreasing $D$ and $\g$.  If Conditions (Z1) and (Z2) are satisfied,  we say that we have mod-Gaussian convergence for the sequence $\left(X_n\right)_{n\in\N}$ with zone of control $[-Dt_n^\g,Dt_n^\g]$ and index of control $(v,w)$. 
\end{definition}

\begin{thm}[Speed of convergence, Theorem 2.16 in \cite{Speed}]\ \label{speed}\\
Assume that conditions (Z1) and (Z2) hold, together with the inequality
 $\g\le \frac{v-1}{2}$. Then,
\begin{equation*}
d_{Kol}\left(\frac{X_n}{\sqrt{t_n}}, \mathcal{N}(0,1)\right)\le C(D, v,K_1)\frac{1}{t_n^{\g+\half}},
\end{equation*}
where $d_{Kol}(\cdot, \cdot)$ is the Kolmogorov distance and 
\begin{equation}
\label{defc}
C(D,v,K_1) = \frac{3}{2 \pi} \left(2^{v-1} \G\left(\frac{v}{2} \right)K_1 + \frac{7}{D} \sqrt{\frac{\pi}{2}}\right).
\end{equation}
\end{thm}

\begin{thm}[Local limit theorem, \cite{DBMN17}]\ \label{llt}\\
 Let $x\in\R$ and $(a,b)$ be a fixed interval, with $a<b$. Let $(X_n)_{n \in \mathbb N}$ be a sequence that converges mod-Gaussian. Assume that 
 conditions (Z1) and (Z2) are satisfied. Then for every $\delta\in\left(0,\frac{1}{2}+ \gamma\right)$, then
 \begin{equation*}
\lim_{n\to\infty}(t_n)^{\d}\P\left(\frac{X_n}{\sqrt{t_n}}-x \in \frac{1}{(t_n)^{\delta}}(a,b)\right)= (b-a)\frac{e^{-\frac{x^2}{2}}}{\sqrt{2\pi}}.
\end{equation*}
In particular, assuming $\g>0$ and taking  $\d=\half$ one obtains
 \begin{equation*}
\lim_{n\to\infty}(t_n)^{\half}\P\left[X_n-x\sqrt{t_n} \in (a,b)\right]=(b-a)\frac{e^{-\frac{x^2}{2}}}{\sqrt{2\pi}}.
\end{equation*}
\end{thm}
We shall see that all the eigenvalues statistics considered satisfy the definition (\ref{def:modGauss}), together with conditions \ref{Z1} and \ref{Z2}.

\section{The eigenvalue density of beta ensembles}\label{subs}
To simplify the writing of most formulas, we adopt the notation:
\begin{align}
\beta' = \betah\,.
\end{align}

\subsection{Notations and previous results}
\label{previous}
We first consider particles living in $\mathbb R$. 
It is known (Section 2.6 in \cite{agz}) that under assumption (\ref{confin}), the empirical  spectral distribution 
\begin{equation}
\label{Ln}L_n = \frac{1}{n} \sum_{k=1}^n \delta_{\lambda_k}\end{equation}
converges in probability towards a compactly supported probability measure $\mu_V$ which is the unique minimizer of
\[\mathcal E_V (\mu) = \int V(x) d\mu(x) - \int\! \int \log |x-y| d\mu(x) d\mu(y)\,.\] If moreover $V$ is real-analytic then the support consists in a finite disjoint union of compact intervals \[{\bf S} = \cup_{h=0}^q [\alpha_h^-, \alpha_h^+]\] and $\mu_V$ has a density of the form
\[\rho_V (x) = \frac{P(x)}{\pi} \prod_{h=0}^q \sqrt{(\alpha_h^+ -x) (x - \alpha_h^-)},\]
where $P$ is analytic in a neighborhood of ${\bf S}$.
It is known that $(L_n)$ satisfies the LDP at scale $n^2$ with good rate function
\begin{equation}
\label{defrate}
\mathcal I(\mu) = \mathcal E_V(\mu) - \mathcal E_V (\mu_V)\,,
\end{equation}
and  that
 \[\lim \frac{1}{n^2} \log \mathcal Z_n^V (\beta) = - \beta' \mathcal E_V(\mu_V)\,,\]
(see \cite{agz}). Furthermore, precise expansions of $\log {\mathcal Z}_n^V(\beta)$ are given by Shcherbina (\cite{sh13}), Borot-Guionnet (\cite{BG13}) and Lebl\'e-Serfaty (\cite{LS15}).

Finally, the two last assumptions are 
\begin{itemize}
\item 
(off-criticality) $P(x) > 0$ for $x \in {\bf S}$
\item (control of large deviations) the effective potential 
\[\mathcal F_V (x) = V(x) - 2 \int \log |x-y| d\mu_V(y)\] achieves its minimum value on ${\bf S}$ only.
\end{itemize}
To establish the AEP property, Bufetov et al. \cite{BMSS13} considered the log-likelihood
\begin{align}\label{L_def}
\mathcal L_{n, \beta}^V (\lambda) :=  \log {\bf P}_{n, \beta}^V (\lambda)
\end{align}
and claimed (see Theorem 2.1)\footnote{\label{a.s.} The convergence is said to be "almost surely", if we consider $\otimes_{n=1}^\infty \mathbb R^n$ as the ambient probability space equipped with the product probability. Since it is not always natural, we can replace it by "convergence in probability".}
\begin{equation*}
\lim_{n \to \infty} \frac{1}{n} \mathcal L_{n, \beta}^V = - E_\beta^V, \end{equation*}
and (see Theorem 2.3)
\begin{equation*}
\label{B2}n^{-1/2} (\mathcal L_n + n E_\beta^V)\Rightarrow \mathcal N(0 ; \sigma_\beta^2(V)),\end{equation*}
where $E_\beta(V)$ and $\sigma_\beta^2(V)$ are suitable constants.

It could be related to a Gibbsian view of the model, writing
\[{\bf P}_{n, \beta}^V = \frac{1}{\mathcal Z_n^V(\beta)} e^{-\beta \mathcal E_n },\]
where 
the energy is
\[\mathcal E_n = - \sum_{1\leq j < k\leq n} \log |\lambda_j -\lambda_k|+ \frac{n}{2}  \sum_{k=1}^n V(\lambda_k).\] 
Popescu  \cite{POP07} considered 
\begin{align*}
E_n = \frac{1}{n} \sum_{k=1}^n V(\l_k) - \frac{2}{n(n-1)} \sum_{1\leq j < k\leq n} \log |\lambda_j -\lambda_k\textcolor{red}{)}|
\end{align*}
which is a slight variant of  $2\mathcal E_n / n^2$. He proved that when $V=H$, 
\begin{enumerate}
\item[(P1)]
 \begin{equation}
\label{P1} \lim_{n \to \infty} n(E_n - a(n)) = 0
\end{equation}
for  $a(n)$ some deterministic function of $n$, and with the convergence
in the sense of the footnote \ref{a.s.};
\item[(P2)]
\begin{equation}
\label{P2}
n^{3/2} (E_n - a(n))\Rightarrow \mathcal N(0; \sigma_{\beta, H}^2)
\end{equation}
for a suitable $\sigma^2_{\beta, H}$;
\item[(P3)]  
$n( E_n - a(n))$ satisfies an LDP at scale $n$ with an explicit rate function.
\end{enumerate}

In Subsection \ref{subs_bdensity} we extend the above results in two directions: a LDP for general potentials and mod-Gaussian convergence in the Hermite, Laguerre cases.

In the same vein, we will consider ensembles whose particles live in $\mathbb T$.
\subsection{Limiting results of the eigenvalue density of $\beta$ ensembles}
\label{subs_bdensity}
In the sequel we will make  a constant use of the function
\[\ell(x) = \log \Gamma(x)\] and we denote 
\begin{equation}
\label{deff0}
f_0(\beta) = \log 2\pi - \ell\left(1 + \beta'\right) + \beta'\log \beta' - \half\left(1 + \beta'\right)\,.
\end{equation}
In the first subsection,  an important role is played by  the Shannon entropy  of a probability distribution $d\m(x) = \mu'(x) dx$ defined by
\[\mathcal S(\mu) = \int \mu'(x) \log \mu'(x) dx\,.\]
When $\mu$ is the equilibrium measure $\mu_V$ associated with a potential $V$ we denote $\mathcal S(V)$, and in the classical cases we denote $\mathcal S(H), \mathcal S(L),$ etc. 
\subsubsection{Real line}
\begin{thm}[LDP]\label{LDPAEPr} $\\$
\begin{enumerate}
\item
The sequence $\left(n^{-1}\mathcal L_{n, \beta}^H\right)_{n \in \N}$ satisfies the LDP at scale $n$ and good rate function
\begin{equation}
\Lambda_{\b, H}^*(x) = \sup_{t> -1} \left(tx - \Lambda_{\b,H} (t) \right),
\end{equation}
where for $t > -1$
\begin{equation}
\label{defLambdaH}
\Lambda_{\beta,H}(t) =
f_0((1+t)\beta) - (1+t) f_0(\beta)\,.
\end{equation}
Moreover, $n^{-1} \mathcal L_{n, \beta}^H$ converges in probability to $-E_\beta^H$ where
\[E_\beta^H := \log 2\pi - \frac{1}{2}  - \beta' + \beta'\Psi\left(1 + \beta'\right) - \ell \left(1+ \beta'\right)\]
where $\Psi$ is the Digamma function.
\item For $V$ satisfying assumptions of Section \ref{previous}, the sequence  $\left(n^{-1}\mathcal{L}_{n, \beta}^V\right)_{n \in \N}$ satisfies the LDP at scale $n$ and good rate function
\begin{equation}
\Lambda_{\b,V}^*(x) = \Lambda^*_{\b,H} (x+ \Delta_{H,V})
\end{equation}
with
\begin{equation}
\label{diffent}
\Delta_{H,V} = \mathcal S(H) - \mathcal S(V)\,.
\end{equation}
Moreover,  
$n^{-1} \mathcal L_{n, \beta}^V$ converges in probability to $-  E_\beta^V $ where

\begin{align}
E_\beta^V = E_\beta^H + \Delta_{H,V}
\,.
\end{align}
\end{enumerate}
\end{thm}
\begin{rmk}
The  value of $E_\b^V$ is proposed in formula (27) in \cite{BMSS13} but there is a mistake
in the expression of $\Delta_{H,V}$ (denoted $S(\rho_V)$ there)  in \cite{BMSS13} (formula $(25)$),   \cite{sh11} (formula $(1.31)$) and \cite{sh13} (formulas $(1.20)$, $(1.21))$. 
Indeed, taking $H =V$ we have to recover $S( \rho_H)= 0$ which is not true, since 
$\mathcal S(H) = \frac{1}{2} - \log  2\pi$
as it is stated in the following proposition.
\end{rmk}
\begin{proposition}
\label{deltas}
\begin{enumerate}
\item When $V(x) = H(x) = x^2/2$,
\begin{align} 
\label{deltaH}
\mathcal S(H) = \half - \log 2\pi\,.
\end{align}
\item When $V(x) = L(x) = \theta x - (\theta -1) \log x$, 
\begin{align}
\label{deltaL}
\mathcal S(L) = 1 - \log 2\pi  + \frac{(\theta -1)}{2} \log (\theta -1) - \frac{(\theta-2)}{2} \log \theta\,.
\end{align}
\item When $V(x) = J(x) = - \kappa_1 \log x - \kappa_2 \log (1-x)$,
\begin{align}
\nonumber
\mathcal S(J) =& - \log (2\pi) +\half\left(h(\kappa_1) + h(\kappa_2) - h(1+ \kappa_1)- h(1 + \kappa_2)\right)\\
\label{deltaJcomplete}
& +\frac{3}{2}\left(h((2+ \kappa_1 + \kappa_2) - h(1+ \kappa_1 + \kappa_2)\right)\\
\end{align}
where $h(x) := x \log x$.
\item When $V(x) = (\dd + 1)\log(1+x^2)$, 
\begin{align}
\label{deltaCy}
\mathcal S(Cy) =& -\log \pi - (1+ \dd) \log (1+ \dd)  + 3 \left(\dd + \half\right) \log \left(\dd + \half\right) - \half \log 2\,.
\end{align}
\end{enumerate}
\end{proposition}
\medskip

Let us present now the results on mod-Gaussian convergence. Set
\begin{equation}
Y_{n, \b}^V = \mathcal L_{n, \b}^V + n E_\b^V\,.
\end{equation}

\begin{thm}[Mod-Gaussian convergence]\ \\
\label{modGAEPr}
The sequence $\left(n^{-1/3} Y_n^{H,\b} \right)_{n \in \N}$ converges mod-Gaussian on the strip $S_1$,  as $n \to \infty$ with the parameters $t_n^\b := n^{1/3} \sigma^2_{\b}$  and the limiting function 
\begin{equation}
\label{AEPpsi}
\psi^{\b}(z) := e^{-A_\beta \frac{z^3}{6}}\,,
\end{equation}
where 
\begin{equation}
\label{kbeta}\sigma^2_{\b} = \beta' - \beta'^2\ell''\left(1+\beta'\right) \ , \ A_\beta =  \beta'^3 \ell'''(1+\beta') + \beta'\,.
\end{equation}
\end{thm}
\begin{conjecture}
The above theorem is still valid, if $H$ is replaced by a potential $V$ satisfying assumption of Section \ref{previous}, with $\sigma_\b^2$ and $A_\b$ as above.
\end{conjecture}
This conjecture comes from an examination of the proof of Theorem \ref{modGAEPr} and from a rough application of formula (\ref{H2V}), connecting the general case with the Hermite case. Actually we are able to claim that this conjecture is true in the Laguerre case. It is also likely that it is true for  the Jacobi and generalized Cauchy cases, since in these cases there is also a Selberg formula. 

\begin{thm}
\label{modGAEPL}
The result of Theorem \ref{modGAEPr} is still valid in the Laguerre case, with $H$ replaced by $L$.  
\end{thm}
\begin{rmk}
For the model of Laguerre potential given by Bufetov et al. (\cite{BMSS13} formula (5)) the normalization is different, but a similar computation 
leads to the same conclusion, the difference is only in the asymptotical expectation.
\end{rmk}

\subsubsection{The circular case}

In the same vein, we have the following result for the circular case.
\begin{thm}
\label{AEPC}
\begin{enumerate}
\item
The sequence $(n^{-1}\mathcal L^C_{n,\b})_{n \in \mathbb N}$ satisfies the LDP at speed $n$ with good rate function
\begin{equation}
\label{rateC}
\Lambda_{\b, C}^* (x) = \Lambda_{\b, H}^* \left(x + \half\right)
\end{equation}
\item
Let
\[Y_{n,\b}^C = \mathcal L_{n,\b}^C + nE_\b^C\,.\]
Then the sequence $\left(n^{-1/3}Y_{n,\b}^C\right)$ 
converges mod-Gaussian on the strip $S_1$,  as $n \to \infty$ with the parameters $t_n^\b := n^{1/3} \sigma^2_{\b}$  and the limiting function $\psi^{\b}$ given by (\ref{AEPpsi}) and (\ref{kbeta}), respectively.
\end{enumerate}
\end{thm}

\begin{rmk}
In the real case, the extension from the Hermite case to the general case is a consequence of Corollary 1.4 in \cite{LS15}. In the circular case, we guess that a similar formula is can be established, which would entail that the LDP holds for general potentials in the circular case. Besides, an explicit computation can be made in the circular Jacobi case, starting from (\ref{beta1+z}) and (\ref{HP}). We let it to the reader.
\end{rmk}
\subsubsection{Consequences}
The following statements hold for $i = H, L, C$
\begin{thm}[Precise moderate deviations] \ \\
For $x > 0$, 
\begin{align*}
&\P\left[\frac{Y_n^{i,\b}}{ n^{2/3}  \sigma^2_{\b}} \ge x\right] =\frac{e^{-\frac{x^2n^{1/3}\sigma^2_\b}{2} }}{x \sigma_\b\sqrt{2\pi n^{1/3}}}\psi^{\b}(x)\left(1+o(1)\right),
\end{align*}
and for $x<0$,
\begin{align*}
&\P\left[\frac{Y_n^{i,\b}}{ n^{2/3}  \sigma^2_{\b}} \le x\right] =\frac{e^{-\frac{x^2n^{1/3}\sigma^2_\b}{2} }}{|x|\sigma_\b\sqrt{2\pi n^{1/3}}}\psi^{\b}(x)\left(1+o(1)\right).
\end{align*}
\end{thm}
\begin{thm}[Extended central limit theorem]\ \\
For $y=o\left(n^{1/6}\right)$, 
\begin{equation}
\P\left[\frac{Y_n^{i,\b}}{\sigma_\b\sqrt{n}}\ge y\right]=\P\left[\mathcal{N}(0,1)\ge y\right]\left(1+o(1)\right).
\end{equation}
\end{thm}

Moreover, it is possible to show that we have a zone of control associated with the mod-Gaussian convergence of the sequences
$\left(n^{-1/3} Y_n^{H, \beta}\right)_{n\in \mathbb N}$.
This, in turn, implies the following additional quantitative estimates.
\noindent\begin{thm}
\label{newthm}
\begin{enumerate}
The following is valid for $i = H,C$.
\item We have, for some constant $K > 0$
\begin{align}
d_{Kol} \left(\frac{Y_n^{i, \beta}}{\sigma_\beta\sqrt n} ; \mathcal N(0,1)\right) \leq K n ^{-1/6}\,.
\end{align}
\item
Let $x\in\R$ and $(a,b)$ be a fixed interval, with $a<b$. Then for every $\d\in\left(0,\frac{3}{4}\right)$,
\begin{align}
\lim_n \left(\sigma^2_\beta n^{1/3}\right)^\d \P\left( \frac{Y_n^{H, \beta}}{\sigma_\beta\sqrt n} -x\in \frac{1}{\left(\sigma^2_\beta n^{1/3}\right)^\d}(a,b) \right) = (b-a)\frac{e^{-\frac{x^2}{2}}}{\sqrt{2\pi}}\,.\end{align}
\item
Let $x\in\R$ and $(a,b)$ be a fixed interval, with $a<b$. Then for every $\d\in\left(0,\frac{3}{2}\right)$,
\begin{align}
\lim_n \left(\sigma^2_\beta n^{1/3}\right)^\d \P\left( \frac{Y_n^{C, \beta}}{\sigma_\beta\sqrt n}-x \in \frac{1}{\left(\sigma^2_\beta n^{1/3}\right)^\d}(a,b) \right) = (b-a)\frac{e^{-\frac{x^2}{2}}}{\sqrt{2\pi}}\,.
\end{align}
\end{enumerate}
\end{thm}

\section{Proofs of Theorems of Section \ref{subs_bdensity}}
\label{proofdens}

\subsection{Preliminaries}
We have to consider, for $\Re z > -1$
\begin{equation}
\label{beta1+z}\log \E \left[e^{ z \mathcal L_{n\beta}^V} \right] = \log \mathcal Z_n^V(\beta(1+z)) - (1+z) \log \mathcal Z_n^V(\beta)\,.\end{equation}
Following Corollary 1.4 in \cite{LS15}, one can see that
\begin{eqnarray}
\label{H2V}\nonumber
\log \mathcal Z_n^V (\beta) &=&\log \mathcal Z_N^H (\beta) - \frac{\b}{2} n^2[\mathcal{E}_V(\mu_V) - \mathcal{E}_H(\mu_H)] +n \left(1-\frac{\b }{2}\right) \Delta_{H,V}\\ 
&& + o((1+\beta)n)\end{eqnarray}
so that, it is necessary to look at the reference Hermite case first.

When $V(x) = x^2/2$,
 the partition function is (\ref{SelbergH}).
Using classical approximations of the Gamma function, 
Bufetov et al. \cite{BMSS13} gives the expansion
\begin{align}\label{7.3}
\log \mathcal Z_n^{H} (\b) =-\frac{3}{4} \beta' n^2 + \beta'  \log n + n f_0(\b) + O(\log n),
\end{align}
where  this $O$ depends on $\beta$.
\subsection{Large deviations}
\subsubsection{Real line : proof of Theorem \ref{LDPAEPr}}
The expansion (\ref{7.3})  gives easily, for $t > -1$
\begin{equation}
\label{6.4}
\lim \frac{1}{n} \log \mathbb E [\exp t\mathcal L_{n,\b}^H]  = f_0 ((1+ t)\beta) - (1+t) f_0(\beta) =: \Lambda_{\b, H}(t)\,. 
\end{equation}
This allows to apply G\"artner-Ellis theorem (\cite{demboz98} Th. 2.3.6). The  function $ \Lambda_{\b, H}$ is steep, i.e. 
\begin{equation}
\lim_{t\downarrow -1} \frac{d\Lambda_{\b, H}(t)}{dt} = -\infty
\end{equation}
 since  $\frac{d\Lambda_{\b, H}(t)}{dt} = \beta' (1+ o(1)) \log (1+t)$ as $t\downarrow -1$.
As a result, we deduce the convergence in probability
\[\frac{1}{n}\mathcal L_{n, \b}^H \rightarrow \frac{d\Lambda_{\b, H}(t)}{dt}\big|_{t=0}  = -\log 2\pi + \frac{1}{2} + \beta' - \beta'\Psi(1+\beta') + \ell (1+\beta') =-  E_\beta^H
\,.\]
 In the general case, using (\ref{H2V}) we have 
\begin{align}
\label{deflambdaVbeta}
\lim \frac{1}{n} \log \mathbb E [\exp t\mathcal L_{n,\b} ^V] = \Lambda_{H,\b}(t)  -t \Delta_{H,V} =: \Lambda_{V, \beta}(t)\,. \end{align}
Its derivative at $0$ is
\[-E_\beta^H - \Delta_{H,V} = -E_\beta^V\,.\]
Consequently, the rate function of the LDP is 
\[\Lambda_{V,\b}^* (x) = \Lambda_{\b, H}^* (x+\Delta_{H,V})\] 
which implies the convergence in probability
\[\frac{1}{n}\mathcal L_{n, \b}^V \rightarrow -E_\beta^V\,.\]

\subsubsection{Unit circle : proof of Theorem \ref{AEPC} 1)}
We start again from the Selberg formula (\ref{SelbergC})
which yields
\begin{align}
\label{circ}
&\log \mathcal Z_n^C(\b) =n\log 2 \pi +\ell\left( 1+n\beta'\right) - n \ell\left( 1 + \beta'\right)\,,
\end{align}
and 
\begin{eqnarray}
\nonumber
\log \mathbb E [\exp t \mathcal L_{n\b}^C] &=& -n t \log 2\pi + \ell (1 + n\b'(1+t)) - (1+t) \ell(1+n\b')\\
&& -n \ell(1 + \b'(1+t)) + n (1+t)\ell(1+\b')
\end{eqnarray}
so that, using Binet,
\begin{eqnarray}
\nonumber
\lim_n \frac{1}{n} \log \mathbb E [\exp t \mathcal L_{n\b}^C]  &=& -t \log 2\pi + \beta'(1+t) \log(1+t) -\ell (1+ \beta'(1+t))\\
\nonumber
&&+ (1+t) \ell(1+ \beta')\\
&=& \Lambda_{\b, H}(t) - \frac{t}{2} \,.
\end{eqnarray}
By duality, we get (\ref{rateC}).

\subsection{Mod-Gaussian convergence}
\subsubsection{Hermite case : proof of Theorem \ref{modGAEPr}}
We use the fundamental relation (\ref{beta1+z}).
The expansion given in Bufetov et al. is not sufficient for us since the dependence on $\b$ of the remaining parts is not controllable. 

We have to start again from the Selberg formula (\ref{SelbergH}) and carefully follow the dependence on $\b$.
This yields
\begin{eqnarray}\nonumber
\log {\mathcal Z}_n^{H}(\b)   = \frac{n}{2} \log 2 \pi + \frac{1}{2}((\beta' - 1) n -\beta' n^2) \log (n\beta') + \sum_{j=1}^n \ell(1+\beta' j) - n \ell(1+\beta')\,.\\
\label{expansion_partition}
\end{eqnarray}
Going to the first Binet formula (\ref{Binet1bis}),   
we have
\begin{eqnarray}\nonumber\ell (1+ \beta' j) &=& \left(\beta'j + \half\right)\log (\beta' j) - \beta' j + \half \log 2\pi + \frac{1}{12\beta' j} - \int_0^\infty \phi(s) e^{-s\beta'j} ds\\
\nonumber
&=&\beta'\ell(1+j) 
 + \frac{1-\beta'}{2} \log (2\pi j) + \left(\beta' j + \half\right) \log \beta'
+ \frac{1}{12j}\left(\frac{1}{\beta'} - \beta'\right)\\
\label{binet}
&&+ \int_0^\infty \phi(s)  \left(\beta' e^{-sj} - e^{-s\beta'j}\right) ds\,.
\end{eqnarray}
Adding up we get
\begin{eqnarray}
\nonumber
\sum_{j=1}^n \ell (1+ \beta' j) &=&\beta' \log G(n+1) + \frac{1+\beta'}{2}\ell (1 +n) 
\\
\nonumber
&&+ \frac{n}{2}\left((n+1) \beta'\log \beta' + (1-\beta')\log 2\pi + \log \beta'\right)
\\
\nonumber
&& + \frac{1}{12}\left(\frac{1}{\beta'}-\beta'\right)\left[\log n + O_n(1) 
\right]\\
&& + \int_0^\infty \phi(s) \left(\frac{\beta'(1 - e^{-ns})}{e^s - 1} - \frac{1-e^{-n\beta' s}}{e^{s\beta'}- 1}\right) ds\,,
\end{eqnarray}
where $O_n(1)$ is bounded and independent of $\beta$.

Using the Binet formula (\ref{Binet1bis}) again and the expansion of Barnes function (\ref{4.184}), we obtain
\begin{eqnarray}
\nonumber
\sum_{j=1}^n \ell (1+ \beta' j) &=& \frac{\beta'}{2}n^2 \log n + \left(\beta' \log \beta'- \frac{3\beta'}{2}\right)\frac{ n^2}{2} + \frac{1+\beta'}{2}  n\log n\\
\nonumber
&& + \left(\beta' \log \beta' + \log \beta'  -1 -\beta' + \log 2\pi\right) \frac{n}{2} \\
\nonumber
&&+ R(\beta) \log n + \beta'[\zeta'(-1) + \varepsilon_{2,n}]\\
\label{expansionell}
&&  + \frac{O_n(1)
}{12}\left(\frac{1}{\beta'}-\beta'\right) + \frac{1+\beta'}{4}\log 2\pi + r_n(\beta),
\end{eqnarray}
where the term $\varepsilon_{2,n}$ denotes the reminder coming from Barnes function expansion (\ref{4.184}),
\[R(\beta) = \frac{1}{4} + \frac{\beta}{24} + \frac{1}{6\beta}\]
and
\begin{equation}r_n(\beta) = \int_0^\infty \phi(s)\left(\frac{\beta'(1 - e^{-ns})}{e^s - 1} - \frac{1-e^{-n\beta' s}}{e^{s\beta'}- 1} - \frac{1+\beta'}{2} e^{-ns}\right) ds\,.\end{equation}
 This yields
\begin{eqnarray}
\nonumber
\log {\mathcal Z}_n^{H}(\b)   &=& -\frac{3\beta'}{4} n^2 +\beta' n\log n +
f_0(\beta) n +R(\beta) \log n\\
\nonumber
&& + \beta'[\zeta'(-1) + \varepsilon_{2,n}] + \frac{O_n(1)
}{12}\left(\frac{1}{\beta'}-\beta'\right) + \frac{1+\beta'}{4}\log 2\pi\\
\label{expansionH}
&&+ r_n(\beta)
\end{eqnarray}
with $f_0$ defined in (\ref{deff0}). 

Now we are ready to use (\ref{beta1+z}). To compute $\log \mathcal Z_n^H (\beta(1+\zeta)) - (1+\zeta) \log \mathcal Z_n^H(\beta)$, with $\Re \zeta > -1,$ we extend (\ref{expansionH}) to the complex domain.
The first two terms of the expansion are canceled.  The contribution of the third term (after scaling by $n$) 
gives $\Lambda_{\beta , H}(s)$ as defined in (\ref{defLambdaH}).
We obtain
\begin{eqnarray}
\label{decomp}
\log \mathbb E \left[ e^{\zeta \mathcal L_{n, \b}^H} \right] = \Lambda_{\b, H}(\zeta) n + B(\zeta)\log n + C_n(\zeta) + D_n(\zeta),
\end{eqnarray}
where 
\begin{eqnarray}
\label{defB}
B(\zeta) &:=& -\zeta\left(\frac{1}{4} + \frac{2+\zeta}{12\beta'(1+\zeta)}\right) ,\\
\nonumber
C_n(\zeta) &:=& - \frac{O_n(1)
}{12\beta'} \frac{\zeta (2+\zeta)}{1+\zeta} - \frac{\zeta}{4}\log 2\pi,\\
\nonumber
D_n(\zeta) &:=& r_n(\beta(1+\zeta)) -(1+\zeta)r_n (\beta).
\end{eqnarray}
Note when $\zeta \rightarrow 0$, $B(\zeta) = O(\zeta)$ and 
\begin{eqnarray}
\label{azeta}\Lambda_{\beta, H} (\zeta) &=& \left(\beta f'_0(\beta) - f_0(\beta)\right)\zeta + \frac{\beta^2 \zeta^2}{2}f_0^{(2)}(\beta) + \frac{\beta^3 \zeta^3}{6}f_0^{(3)}(\beta) + O(\zeta^4)
\end{eqnarray}
and  $C_n(\zeta) = O(\zeta)$. 

Let us study the remaining term. We have
\begin{eqnarray}
\nonumber
D_n(\zeta) &=& r_n((1+\zeta)\beta) - (1+\zeta) r_n(\beta)\\
\label{defDn}
&=&\int_0^\infty \phi(s) \left(\frac{1 - e^{-n\beta's}}{e^{s\beta'} - 1}(1+\zeta) - \frac{1 - e^{-n\b'(1+\zeta)s}}{e^{s\b'(1+\zeta)} -1} +\frac{\zeta}{2} e^{-ns}\right) ds\,. 
\end{eqnarray}
Since $\phi(s) \leq 1/12,$ the last contribution yields
\begin{eqnarray*}
\frac{|\zeta|}{2} \int_0^\infty \phi(s) e^{-ns} ds \leq \frac{|\zeta|}{\textcolor{red}{24}n}\,.
\end{eqnarray*}
In the same way, owing to (\ref{720}),  
\begin{eqnarray*}
|\zeta|\int_0^\infty \phi(s) \frac{1- e^{-n\beta's}}{e^{s\beta'} - 1} ds \leq |\zeta|\int_0^\infty  \frac{\phi(s)}{e^{s\beta'} - 1} ds \leq K|\zeta|\,.
\end{eqnarray*}
Now
\begin{eqnarray}
\nonumber
\frac{1 - e^{-n\beta's}}{e^{s\beta'} - 1} - \frac{1 - e^{-n\b'(1+\zeta)s}}{e^{s\b'(1+\zeta)} -1} = D^1_n(s) - D^2_n(s)\,,
\end{eqnarray}
where
\begin{eqnarray}
\nonumber
D^1_n(s) =\frac{e^{-n\b's}(e^{-n \b' \zeta s} -1)}{e^{s\b'(1+\zeta)} -1}
\end{eqnarray}
and
\begin{eqnarray}
\nonumber
D^2_n(s) = \frac{(1- e^{-n\b's})e^{\b's} (e^{\b'\zeta s} -1)}
{\left(e^{s\b'(1+\zeta)} -1\right)\left(e^{s\b'}- 1\right)}\,.
\end{eqnarray}
 We use two elementary bounds for $z \in \mathbb C$:
\begin{align}
\label{elemin}
|e^z -1| &\geq |e^{\Re z}-1|\\
\label{elemaj}
|e^z -1 |& \leq |z| \max(e^{\Re z},1)\,.
\end{align}
Now, if we assume
\begin{equation}
\label{bigbound}
|\zeta| \leq a < 1\,,
\end{equation}
then the above bounds give, for every $s \geq 0$
\begin{eqnarray}
|e^{s\b'(1+\zeta)} -1| \geq e^{s\b'(1-a)} -1
\end{eqnarray}
\begin{equation}
e^{-n\b's}|(e^{-n \b' \zeta s} -1)|\leq e^{-n\b's} n\b's|\zeta| e^{n\b'as}\,.
\end{equation}
Integrating and using
\[\max_{x >0} xe^{-x} = e^{-1}\,,\]
we get
\begin{equation}
\int_0^\infty |D^1_n(s)|  \phi(s) ds \leq |\zeta|\frac{e^{-1}}{1-a}\int \frac{\phi(s)}{e^{s\b'(1-a)} -1} ds = K_a^1 |\zeta|\,.
\end{equation}
In the same vein
\begin{equation}
\int_0^\infty |D_n^2(s)| \phi(s) ds \leq |\zeta| \int \frac{s e^{\b'(1+a)s}\phi(s) }{(e^{s\b'(1-a)} -1)^2} ds = K_a^2 |\zeta|\,,
\end{equation}
which allows to get, for $|\zeta| \leq a  <1$
\[|D_n(\zeta)|\leq \left(K+ K_a^1 + K_a^2 + \frac{1}{12 n}\right) |\zeta| \,.\]
Besides
\begin{eqnarray}
\label{betaf}
\nonumber
\beta f'_0 (\beta) - f_0(\beta) &=& -\log 2\pi + \half + \b' - \b' \Psi(1+\b') + \ell (1+\b') = - E_\beta^H\\
\nonumber
\beta^2f_0^{(2)}(\beta) &=&  
\beta' - \left(\beta'\right)^2\Psi'\left( 1 + \beta'\right) = \sigma^2_\beta\\
\beta^3f_0^{(3)}(\beta) &=&
 -\beta'  - \left(\beta'\right)^3\Psi''\left( 1 + \beta'\right)= -A_\beta\,.
\end{eqnarray}
Gathering, we get, for $|z| \leq a n^{1/3}$ 
\begin{eqnarray}
\nonumber
\log \mathbb E \left[e^{ zn^{-1/3}\mathcal L_{n,\beta}^H}\right]&=& 
-E_\beta^H n^{2/3} z+\frac{\sigma_\beta^2 z^2}{2} n^{1/3}-A_\beta \frac{z^3}{6}\\
\label{finalH}
&&+ O\left(\frac{|z|^4}{n^{1/3}}\right)+ O \left(\frac{|z| \log n }{ n^{1/3}} \right)  + \frac{K_a|z|}{n^{1/3}},
\end{eqnarray}
which ends the proof.

\subsubsection{Laguerre case : proof of Theorem \ref{modGAEPL}}

We use a careful analysis of the partition function (\ref {SelbergL}). To make easier the calculations, and in particular to use (\ref{expansionell})  we will assume that $n\theta =: q \in \mathbb N$. In that case we may write 
\begin{align}
\mathcal Z_n^L(\beta) = 
\frac{ \left(n\beta' \theta\right)^{-n(\theta n \beta' + 1 - \beta')}\Gamma(1 + \beta'(\theta -1)n)} { \Gamma(1 + \beta')^n\Gamma(1 + \beta'\theta n)}\times 
\frac{\displaystyle \prod_{j=1}^n \Gamma(1+ \beta' j)\prod_{k=1}^{\theta n}\Gamma(1+ \beta' j)}{\displaystyle \prod_{k=1}^{\theta n -n}\Gamma(1+ \beta' j)},
\end{align}
so that
\begin{align}
\log \mathcal Z_n^L(\beta)  = A_n + B_n 
\end{align}
with
\begin{align}
\nonumber
A_n :=& -n(\theta n \beta' + 1 - \beta') \log (n\beta' \theta) - n \ell(1+\beta') - \ell (1 + \beta'\theta n) + \ell(1 + \beta'(\theta -1)n)
\end{align}
and
\begin{align}
B_n := \sum_{j=1}^n \ell (1+ \beta'j) + \sum_{j=1}^{\theta n} \ell (1+\beta'j) - \sum_{j=1}^{(\theta -1)n} \ell(1+\beta'j)\,.
\end{align}
Using the first Binet formula (\ref{Binet1bis}), we derive
\begin{align}
\nonumber
A_n =& - \theta \beta' n^2 \log n - \theta \beta' \log (\beta'\theta) n^2  - n \log n\\
\nonumber
& + \left(-\log \beta' - (1-\beta') \log \theta- \beta'(\theta \log \theta - (\theta -1) \log (\theta -1))
 + \beta' - \ell(1 + \beta')\right)n
\\
\nonumber
&- \half \log (\theta/ (\theta-1)) \\
\nonumber
&+\frac{1}{12\beta'\theta (\theta -1)n}\\
\label{firstest}
&-\int_0^\infty\phi(s)(e^{n\beta' s} -1) e^{-ns\beta'\theta} ds\,.
\end{align}
To ease the reading of the computations of the second estimation, 
let us first notice some elementary facts. If we set, for $u : \mathbb R \rightarrow \mathbb R$
\[T_\theta (u(n)) := u(n) + u(n\theta) - u(n(\theta -1))\,,\]
we have 
\begin{align*}
T_\theta(n^2 \log n) &= 2\theta n^2 \log n +\left(\theta^2\log \theta - (\theta-1)^2 \log (\theta - 1)\right)  n^2,\\
T_\theta(n^2) &= 2 \theta n^2,\\
T_\theta (n \log n) &= 2n \log n + \left(\theta \log \theta - (\theta -1) \log (\theta -1)\right) n ,\\
T_\theta(n) &= 2n,\\
T_\theta (\log n) &= \log n + \log (\theta/(\theta-1)),\\
T_\theta(1) &= 1\,.
\end{align*}
We get from (\ref{expansionell}) 
\begin{align}
\nonumber
B_n
=&\beta' \theta n^2 \log n + \frac{\beta'}{2}\left(\theta^2\log \theta - (\theta-1)^2 \log (\theta - 1)\right)n^2\\
\nonumber
&+\left(\beta' \log \beta' - \frac{3\beta'}{2}\right)\theta n^2\\
\nonumber
&+ (1+\beta')n \log n+ \frac{1+\beta'}{2}  (\theta \log \theta - (\theta -1) \log (\theta -1)) n\\
\nonumber 
&+ (\beta'\log \beta' + \log \beta' -1 - \beta' + \log 2\pi) n\\
\nonumber
&+ R(\beta) \log n + R(\beta) \log (\theta/(\theta-1))\\
\nonumber
&+ \beta'[ \zeta'(-1) + T_\theta (\varepsilon_{2,n}) ] + \frac{T_\theta (O_n(1))}{12}\left(\frac{1}{\beta'} - \beta'\right) + \frac{1+\beta'}{4} \log 2\pi\\
\label{secondest}
&+T_\theta (r_n(\beta))\,.
\end{align}
Joining $A_n$ and $B_n$ gives
\begin{align}
\nonumber
\log \mathcal Z_n^L(\beta)  
&= \frac{\beta'}{2}\left( \theta^2 \log \theta - (\theta-1)^2 \log (\theta - 1) -2 \theta \log \theta - 3\theta\right) n^2\\
\nonumber
&+\beta'  n\log n\\
\nonumber
& + \left(f_0(\beta) + \frac{1 -\beta'}{2} \left[\theta \log \theta - (\theta -1) \log (\theta -1) -1 - 2 \log \theta\right]  \right)n\\
\nonumber
&+ R(\beta) \log n\\
\nonumber
& + \left(R(\beta) - \half\right) \log (\theta/(\theta-1))\\
\nonumber
&+ \beta'[ \zeta'(-1) + T_\theta (\varepsilon_{2,n}) ] + \frac{T_\theta (O_n(1))}{12}\left(\frac{1}{\beta'} - \beta'\right) + \frac{1+\beta'}{4} \log 2\pi\\
\label{globalest}
& + r^L_n (\beta)
\,,
\end{align}
where
\begin{equation}
\label{rL} r^L_n (\beta) =T_\theta (r_n(\beta))+\frac{1}{12\beta'\theta (\theta -1)n}
-\int_0^\infty\phi(s)(e^{n\beta' s} -1) e^{-ns\beta'\theta} ds\,.\end{equation}
Now, as in the Hermite case, we compute $\log \mathcal Z_n^L((1+\zeta) \beta) - (1+\zeta) \log\mathcal Z_n^L(\beta) $ and observe that the first two terms in the expansion (\ref{globalest}) canceled, and the contribution of the coefficient of $n$ is exactly $\Lambda_{\beta, L}(\zeta)$ as given in (\ref{deflambdaVbeta}),  (\ref{diffent}) and (\ref{deltaL})
so that
\begin{align*}
\log \mathbb E \exp(\zeta \mathcal L^L_{n, \beta})&=
\log \mathcal Z_n^L((1+\zeta) \beta') - (1+\zeta) \log\mathcal Z_n^L(\beta') \\&=  \Lambda_{\beta, L}(\zeta) 
n + B(\zeta) \log n + C_n^L (\zeta) + D_n^L (\zeta)
\end{align*} 
where 
$B$ is in (\ref{decomp}) and
\begin{align*}
C_n^L(\zeta) &= \left[B(\zeta) + \frac{\zeta}{2}\right] \log \frac{\theta}{\theta - 1} - \frac{T_\theta(O_n(1))}{12\beta'} \zeta \frac{2+\zeta}{1+\zeta} - \frac{\zeta}{4}\log 2\pi\\
D_n^L (\zeta) &=r_n^L(\beta(1+\zeta)) - (1+ \zeta) r_n^L (\beta)\\
&=T_\theta(D_n(\zeta))- \frac{1}{12\theta(\theta-1)\beta'n} \frac{2\zeta + \zeta^2}{1+\zeta}\\
&+\int_0^\infty \phi(s) \left[e^{-n\theta \beta'(1+\zeta)s}- e^{-n\theta\beta' s}\right]\!\ ds\\
&-\int_0^\infty \phi(s) \left[e^{-n(\theta -1)\beta'(1+\zeta)s}- e^{-n(\theta -1) \beta' s}\right]\!\ ds\\
&+\zeta\int_0^\infty \phi(s) \left[e^{-n(\theta -1)\beta's}- e^{-n\theta \beta' s}\right]\!\ ds\,.
\end{align*}
Arguing as in the previous section, we get
\[C_n^L(\zeta) = O(\zeta)\ , \ D_n^L(\zeta) = O(\zeta)\,.\]
Hence, when $\zeta= zn^{-1/3}$, and using (\ref{azeta})
we conclude
\begin{align*}
\log \mathbb E \left[e^{ zn^{-1/3}\mathcal L_{n,\beta}^L}\right]= -E_\beta^L n^{2/3} z + \frac{\sigma^2_\beta z^2}{2}n^{1/3} - A_\beta \frac{z^3}{6} +  \mbox{a remainder term as in (\ref{finalH})}.
\end{align*}
We skip the end of the computation since it is similar to the Hermite case.\\

\subsubsection{Proof of Proposition \ref{deltas}}
When a measure has a density of the form
\[\rho_V(x) = R(x) \sqrt{((x-a_V)(b_V -x)} \ \ (x \in [a_V, b_V))\]
the function $\log \rho(x)$ is a combination of $\log$ and  the computation of 
$\mathcal S(V)$ lies on the knowledge of the logarithmic potential
\begin{align}\label{logpot0}
U^V(y) = - \int_{\bf S} \log |y-x|\!\ \rho_V(x)\,.\end{align}

\noindent 1. In the Hermite case, we have from \cite{forrester2012spectral} (Prop. 1.1) after scaling 
\begin{align}
\label{ForH}
-U^H(y) = \log 2 - \half + \frac{y^2}{4} - \frac{|y|}{4} \sqrt{y^2 -4} - \log \left[|y| - \sqrt{y^2 -4}\right] \ (|y| \geq 2)\,,
\end{align}
so that $-U^H(\pm 2) = \half$ and then 
\[\mathcal S(H) = |\log 2\pi -\half\left(U^H(2) + U^H(-2)\right) = \half - \log 2\pi\,.\]
\medskip

\noindent 2. In the Laguerre case,  we have
\begin{align}
\label{totaldeltaL}\mathcal S(L) = \log \theta - \log 2 \pi 
- \frac{1}{2} U^L (c_+) - \half U_L(c_-)
 + U_L(0)\,.
\end{align}
It is known from Hiai-Petz \cite{hiai} (formula (5.3.11) after scaling), that
\begin{align}
\label{logLag}
-U^L (x) = \half\left(\theta x  - (\theta - 1) \log x  -  \log \theta   -( 1 + \theta)\right)
 \end{align}
 for $x\in [c_- , c_+]$ . Moreover \begin{align}
\label{latter}
-U^L (0) = -(\theta - 1) \log (\theta - 1) + (\theta -1) \log \theta -1 \,,\end{align}
(see formula (4.6) in \cite{R07} and references therein).
Now, since $\sqrt{c_+c_-} = ((\theta-1)/\theta$,
we get
\begin{align}
\mathcal S(L) =  1  - \log 2\pi  + \frac{(\theta -1)}{2} \log (\theta -1) - \frac{(\theta-2)}{2} \log \theta \,.
\end{align}
Notice that for the computations of $U^L(0)$ and $U^L(c_\pm)$ we could also have started from Prop.  3.1 in \cite{forrester2012large}.
\medskip

\noindent 3. In the Jacobi case, we have
\begin{align}
\label{defdeltaJ}\mathcal S(J) = \log C_{a_-, a_+} - \log 2\pi  + U^J(0) + U^J(1) -\half\left(U^J(a_+) + U^J(a_-)\right)\,.\end{align}
Actually, Forrester in \cite{forrester2012large} formula (3.27) gives the expression of the logarithmic potential.
It is
\begin{align}
\nonumber
-U^J(x) = &\log \half\left|x- \frac{a_+ + a_-}{2} +u\right| - \kappa_1 \log \frac{2|\sqrt{a_+a_-} + x-u|}{\left(\sqrt{a_+} + \sqrt{a_-}\right)^2}\\
\label{ForJ+}
&- \kappa_2 \log \frac{2|\sqrt{(1-a_+)(1-a_-)} + 1 -x +u|}{\left(\sqrt{1-a_+} + \sqrt{1-a_-}\right)^2}\,,
\end{align}
where $u = \pm \sqrt{(x-a_-)(x-a_+)}$ where $\pm$ is $+$ if $x\in [a_+ , 1]$ and $-$ if $x\in[0, a_-]$.
Actually, there are three system of parameters for this kind of study: $(\kappa_1, \kappa_2)$ , $(a_-, a_+)$ and $(\sigma_-, \sigma_+)$ (see 
 \cite{R07}) given by
\[\sigma_\pm = \half \left[1 + \sqrt{a_+a_-} \pm \sqrt{(1-a_+)(1-a_-) }\right]\,.\]
The return to $(\kappa_1, \kappa_2)$ is 
\begin{equation}
\label{comeback}\sigma_+ = \frac{1+ \kappa_1 + \kappa_2}{2 + \kappa_1 + \kappa_2} \ , \ \sigma_- = \frac{1+ \kappa_1}{2 + \kappa_1 + \kappa_2}\,.\end{equation}
In the system $(\sigma_-, \sigma_+)$, using the notation $h(x) = x \log x$, we then have
\begin{align*}
U^J (0)  &= \frac{1}{1- \sigma_+} \left[ h(\sigma_+ + \sigma_- -1) -h(\sigma_+) - h(\sigma_-) \right]\\
U^J (1)  &= \frac{1}{1- \sigma_+} \left[ h(\sigma_+ - \sigma_-)  - h(\sigma_+) -h(1-\sigma_-)\right]\,,
\end{align*}
and 
\begin{align*}
-\half (U^J(a_+) + U^J(a_-) ) &= 
 \frac{1}{2(1- \sigma_+)} \left[h(\sigma_+) +h(\sigma_-) + h(1-\sigma_+)+ h(1-\sigma_-)\right]\\
&-  \frac{1}{2(1- \sigma_+)} \left[ h(\sigma_+ - \sigma_-) + h(\sigma_+ + \sigma_- -1)\right]
\end{align*}
Gathering the different terms and using $C_{a_-, a_+} = (1 -\sigma_+)^{-1}$  we get
\begin{align*}
\mathcal S(J) =& - \log (2\pi) - \frac{3}{2(1- \sigma_+)} h(\sigma_+) - \half\log (1- \sigma_+)  \\
&+\frac{1}{2(1-\sigma_+)}\left[ h(\sigma_+ + \sigma_- -1)
 +  h(\sigma_+ - \sigma_-) -h(\sigma_-) -  h(1-   \sigma_-)\right]\,,
\end{align*}
which coming back to $(\kappa_1, \kappa_2)$ via (\ref{comeback}) gives (\ref{deltaJcomplete}) .
\medskip

\noindent 4. In the generalized Cauchy case, we have
\begin{align}\mathcal S(C_y) =& \log \dd - \log \pi + U^{Cy}(i) + U^{Cy}(-i) -\half U^{Cy}(m) - \half U^{Cy}(-m)
\end{align}
To compute 
\[-U^{Cy} (z) = \int_{-\mm}^\mm \log |z-y| d\mu^{Cy} (y)\]
we consider  derivative i
\[G(z) = \frac{(\dd +1)z -\dd \sqrt{z^2 - \mm^2}}{1+z^2}\,.\]
Setting
\[\sqrt{z^2 -\mm^2} = zt \ , \ zdz =  \frac{\mm^2 t}{(1-t^2)^2} dt\]
we get
\begin{align*}G(z)\!\ dz 
=& \frac{(2\dd +1) t}{(1-t^2)\left(\dd t + (1+\dd)\right)}\!\ dt \\
=&- \frac{dt}{2(t-1)} -\frac{(2\dd + 1) dt}{2 (1+t)}+  \frac{(1+\dd) dt}{t+ \frac{(1+\dd)}{\dd}}
\end{align*}
hence
\[-U^{Cy} (z)   = F(t) + C\,\]
where
\[F(t) = -\half \log |1-t| - \frac{2\dd +1}{2}\log (1+t) + (\dd + 1) \log (\dd +1 + \dd t)\]																
and $C$ is such that
\[\lim_{|z| \rightarrow \infty} U^{Cy}(z) +\log |z| = 0\]
or, in other words
\[\lim_{|t| \rightarrow 1} F(t)  -\log |z(t)| = -C\,.\]
Since
\[\log |z(t)| = \half \log (2\dd +1)  - \half \log|1-t| - \half \log (1+t) - \log \dd\]
we get
\[C = -\frac{2\dd +1)}{2} \log(2\dd +1) +\dd  \log 2 - \log \dd\,.\]
In particular for $z= \pm \mm$ we have $t=0$ so that 
\begin{align}
\half U^{Cy}(\mm) + \half U^{Cy}(-\mm) = - \dd \log 2 +  \log d
\end{align}
Taking successively $z= i$ and $z= -i$ and adding up, we get:
\begin{align*}
\int_{-\mm}^\mm \log (1+x^2) d\mu^{Cy} (x) &=  2 (1+2\dd) \log 2 +  2(1+\dd) \log (1+\dd) \\
 & + 2\dd \log \dd - 2(1+2\dd) \log (1+2\dd) 
\end{align*}
and we conclude
\begin{align}\mathcal S(C_y) =& -\log \pi - (1+ \dd) \log (1+ \dd)  + 3 \left(\dd + \half\right) \log \left(\dd + \half\right) - \half \log 2\,.
\end{align}

\begin{rmk}
In the very specific cases of the semi-circle distribution $\rho_H$ or the arcsine distribution, which corresponds to $\kappa_1 = \kappa_2 = 0$ in the Jacobi case, there is another way to compute the entropy via the classical formula
\[\int \rho(x) \log \rho(x)\!\ dx = \left(\frac{d}{ds}\right)_{s=1} \int \rho(x)^s \!\ dx\,\]
and some elementary properties of the Gamma function. 
\end{rmk}
\subsubsection{Circular case : proof of Theorem \ref{AEPC} 2)}
From (\ref{binet}) and (\ref{circ}) we have
\begin{eqnarray}
\nonumber
\log Z_n^C(\b) &=& \b' \ell(1+n) + n \left(\log 2\pi - \ell(1+\b')
 + \b'\log \b'\right)\\
 \nonumber
&& + \frac{1-\b'}{2} \log n + \frac{1-\b'}{2}\log (2\pi)+\half\log(\b')\\
&& + \frac{1}{12n\b'}-\frac{\b'}{12n}+ \int_0^\infty \phi(s)  \left(\b' e^{-sn }- e^{-s\b'n}\right) ds\,,
\end{eqnarray}
hence
\begin{eqnarray}
\nonumber
\log \mathbb E \exp (\zeta \mathcal L_{n, \b}^C) &=& n \left(\Lambda_{\b, H}(\zeta) - \frac{\zeta}{2}\right)\\
\nonumber
&& -\frac{\zeta}{2}\log n-\frac{\zeta}{2}\log (2\pi \b') + \half \log (1+\zeta)-\frac{\zeta(\zeta+2)}{12\b'n(1+\zeta)}\\
\nonumber
&&
 - \int_0^\infty \phi(s) \left(e^{-s\b'(1+\zeta)n} - (1+\zeta)e^{-s\b'n}\right)ds\,.\\
\end{eqnarray}
Let us call the last two terms $F_n$ and $G_n$, respectively.

The main part is now
\[\left(\Lambda_{\b, H}(\zeta) - \frac{\zeta}{2}\right) n\]
which tells us that everything is similar to the Hermite case with 
\[E_\b^C = E_\b^H + \frac{1}{2}\,.\]
The remaining part has the following behavior:
\[ F_n = O(\zeta \log n)\]
and
\begin{eqnarray*}
|G_n| \leq   |\zeta|\int_0^\infty\phi(s) e^{-s\b'n} ds + \int_0^\infty \phi(s)e^{-sb'n}\left|e^{-s\b'zn} -1\right| ds
\end{eqnarray*}
and a straightforward computation leads for $|\zeta| \leq a < 1$ to
\[|G_n| \leq \frac{|\zeta|}{12\b'n} + \frac{a|\zeta|}{12\b'n (1-a)^2}\,,\]
and we conclude as in previous sections.

\subsection{Proof of Theorem \ref{newthm}}

We will prove that the parameters of the zone of control satisfy $\gamma= 1/4$ (resp.  $\gamma= 1$), $v=1$ and $w=2$, or more precisely that
\begin{align}
|\left[\psi_n (i\xi) - 1\right|\leq K_1 |\xi| e^{K_2|\xi|^3}
\end{align}
as soon as $\xi \in [-t_n^\gamma, t_n^\gamma]$. Of course for the first statement, we replace $\gamma$ by $0$.

In the Hermite case, we have (see (\ref{finalH}))
\begin{align}
\psi_n(i\xi) = \exp \left(-\frac{iA_\beta \xi^3}{6} +u_n\right)
\end{align}
with
\begin{align}
u_n(\xi) =  O(|\xi|^4 n^{-1/3})+ O ( |\xi| (\log n) n^{-1/3})  + K_a|\xi|n^{-1/3}\,.
\end{align}
so that, for some constants $k_1, k_2 > 0$
\begin{align}
|u_n(\xi)|&\leq k_1 (|\xi| + |\xi|^4)\ \ \  &(z\in \mathbb R)\\
|u_n(\xi)| &\leq k_2 \ \ \ &(|\xi| \leq n^{1/12})
\end{align}
Using (\ref{elemaj}) we have

\begin{align}\left|\psi_n (i\xi) - 1\right|&\leq
\left(\frac{|\xi|^3A_\beta}{6} + |u_n|\right) e^{u_n}\\
&\leq  
|\xi|\left(\frac{|\xi|^2A_\beta}{6} +  k_1 +k_1 |\xi|^3\right) e^{k_2}
\end{align}
Now, there exits $k_3 > 0$ such that, for every $\xi \in \mathbb R$
\[\frac{|\xi|^2A_\beta}{6} +  k_1 +k_1 |\xi|^3 \leq k_3 e^{|\xi|^2}\,.\]

In the circular case, we can perform the same estimations, starting from 
\[\psi_n (i\xi) = \exp \left(-\frac{iA_\beta \xi^3}{6} +u_n\right)\]
with
\begin{align}
u_n(\xi) =  O(|\xi| n^{-1/3} \log n)\,,
\end{align}
as soon as $|\xi| \leq a n^{1/3}$. 
This yields the result with $\gamma = 1$, $v=1$ and $w=2$.

\section{Appendix : Some properties of the Gamma and the Barnes G functions}

\label{appendix}
The Gamma function 
defined by
\[\Gamma (z) = \int_0^\infty e^{-t} t^{z-1} dt \ , \ (\Re\!\ z > 0)\,.\]
Its extension as a meroporphic function on $\mathbb C \setminus \{-1, -2, \dots\}$ has the representation
\begin{equation}\label{gamma}
\G\left( 1+z\right) = \prod_{j=1}^\infty \Big\{ \left( 1+\frac{z}{j}\right)^{-1} \left( 1+\frac{1}{j}\right)^z \Big\}\,.
\end{equation} 
The first Binet's representation of the logarithm of the Gamma function 
is given by  (see for instance \cite{EMOT} p.21 or \cite{WW96} p.242)
\begin{align}
\label{Binet1}
\ell (z) = \left( z - \half\right) \log z - z + 1 + \int_0^\infty \varphi(s) \left[ e^{-sz} - e^{-s}\right]ds,  \quad \Re (z) >0.
\end{align}
where the function $\varphi$ is defined as
\begin{equation*}
\varphi(s) = \left[\half - \frac{1}{s} + \frac{1}{e^s-1} \right] \frac{1}{s}
\end{equation*}
and satisfies, for every $s\ge 0$
\begin{equation}
\label{12}
0<\varphi(s)\le \varphi(0)=\frac{1}{12}.
\end{equation}
Moreover, if we set
\[\phi(s) = \frac{1}{12}-\varphi(s) \ , \quad s >0\]
we have that  $0 \leq  \phi(s) < 1/12$ for $s > 0$ and
\begin{equation}\label{720}\lim_{x \rightarrow 0} \phi(x)/x^2 = 1/720\,,
\end{equation}
which leads to another version of (\ref{Binet1}), that is
\begin{equation}
\label{Binet1bis}
\ell (1+z) = \left(z + \half\right)\log z -z + \half \log 2\pi + \frac{1}{12z} - \int_0^\infty \phi(s) e^{-sz} ds\,.
\end{equation}
The second Binet's formula (see for instance \cite{EMOT} p.22 or \cite{WW96} p.245) is
\begin{align}
\label{Binet2}
\ell (z) = \left(z-\half \right) \log z - z + \frac{1}{2} \log (2 \pi) + 2 \int_0^\infty \frac{ \arctan\left(\frac{s}{z} \right)}{e^{2 \pi s} - 1} ds, \quad \Re (z) >0,\end{align}
where, for complex $\zeta$, $\arctan\zeta$ is defined as
\[\arctan \zeta := \int_0^\zeta \frac{dt}{1+t^2}\]
with integration along a straight line.
Both have as a byproduct the classical Stirling-like formula (see Olver \cite{FR97} p.293 or \cite{WW96} p.243)
\begin{align}
\label{Stirling}
\ell (1+z) = \left( z + \half \right) \log z - z + \half \log 2 \pi + \frac{1}{12 z} +  r(z), 
\end{align}
where,
for every $\delta > 0$, 
\[\sup_{|\arg z| \leq \pi-\delta}|z|^3 |r(z)| = R_\delta < \infty\,.\]
The derivative of the logarithm of the Gamma function 
\begin{align}
\Psi(x) = \ell' (x)
\end{align}

 satisfies
\begin{align}
\nonumber\Psi(1+x) &= \frac{1}{x} + \Psi(x)\\
\label{Psig}
\Psi(1) = -\gamma \ &, \ \Psi(1/2)= -2\log 2 - \gamma
\end{align}
where 
$\gamma$ is the Euler constant.

The Barnes $G$-function is the entire solution of the functional equation
\begin{align}
\label{finiteG}
G(z+1) = G(z) \Gamma (z)\,,
\end{align}
it may be represented as a product
\begin{align}
\label{expG}
G \left( 1+z \right) = \left( 2 \pi \right)^{z/2} e^{-(z+1)z/2} \prod_{j=1}^\infty \Big\{ \left( 1+\frac{z}{j}\right)^j \left(1 + \frac{1}{j}\right)^{z^2/2} e^{-z}\Big\}\,.
\end{align}
Its derivative satisfies (see \cite{WW96} p.258)
\begin{align}\label{der-Barnes}
\frac{G'(z)}{G(z)} = (z-1) \Psi(z) - z + \half \log 2\pi + \half\,,
\end{align}
Its logarithm is related to the $\ell$ function as follows (see \cite{Barnes})
\begin{align}
\label{Barnes}
\log G(z+1) = \frac{z(1-z)}{2} + \frac{z}{2} \log 2 \pi + z \ell(z) - \int_0^z \ell(x) dx 
, \quad \Re (z) >0.
\end{align}
Gathering the expansions (\ref{Barnes}) and (\ref{Binet2}), we obtain the following integral representation for all $z\in\C$ with $\Re (z)>0$,
\begin{align*}
\log G(z+1)  
&= \frac{z}{2} - \frac{3z^2}{2} + z \log 2 \pi + z \left(z-\half \right) \log z  -  \frac{z}{2} \left( 1-\frac{z}{2} + (z-1) \log z \right) \\
&+ \frac{z^2}{2} - \frac{z}{2} \log (2 \pi) - 2 \int_0^\infty \frac{  \frac{1}{2} \log \left(z^2+s^2 \right) - s \log s}{e^{2 \pi s} - 1} ds ,
\end{align*}
which gives 
\begin{equation*}
\log G(z+1) = \frac{z^2}{2} \log z -\frac{3}{4}z^2 + \frac{z}{2} \log 2 \pi   - \int_0^\infty   \log \left(1+ z^2s^{-2} \right) \frac{s\!\ ds }{e^{2 \pi s} - 1} .
\end{equation*}
Moreover, there is an asymptotic formula (\cite{F10} 4.184) when $|z| \rightarrow \infty$
\begin{equation}
\label{4.184}
\log G(z+1) =  \frac{z^2}{2}\log z - \frac{3}{4}z^2
 + \frac{z}{2}\log 2\pi - \frac{1}{12}\log z + \zeta'(-1) + o(1) \end{equation}
 where $\zeta$ is the Riemann function.

\bibliographystyle{plain}
\bibliography{Bib_AEP}

\end{document}